# Exploring criteria for designing novel waterbomb tessellations using triangular convex polygons

*S. Deshmukh[1], M. Assis[2]*

**Abstract:** Waterbomb style tessellations have been explored in the past by artists such as Ronald D. Resch, Benjamin Parker and Mitya Miller. Generalised waterbomb tessellations are still underexplored in origami design. We have explored various sets of criteria for generalising waterbomb tessellations in order to enumerate valid patterns. We only consider triangular waterbomb tessellations, other polygons will be explored in future papers. In our search we have uncovered some new waterbomb tessellations, which could offer new uses in representational and geometric origami design. We conclude by discussing foldability properties and possible generalisations.

## 1   Introduction

The waterbomb family of origami tessellations is a relatively unexplored area of origami design. The name 'waterbomb' is derived by reference to the traditional waterbomb base which gets its name from its use in the traditional origami model of a waterbomb, see Figure 1. The waterbomb base is a simple collapse of a square such that the square edges become co-planar and there results four connected flaps. We can generalise this concept to other polygons. We call a generalised waterbomb any polygon with creases such that after the polygon is folded, all of the polygon edges are co-planar. In other words, if the folded generalised waterbomb is dipped into a shallow vat of paint, only the polygon edges become coloured. When such polygons can be tiled, either with or without gaps, with constraints discussed below, the resultant crease patterns (CPs) are called waterbomb tessellations, and the gaps in the tiling are called the floors of the tessellation, see Figure 2.

[1] School of Design, RMIT University, Melbourne VIC Australia, corresponding author.
   Email: sukanya.deshmukh@rmit.edu.au
[2] University of Melbourne, Parkville VIC Australia

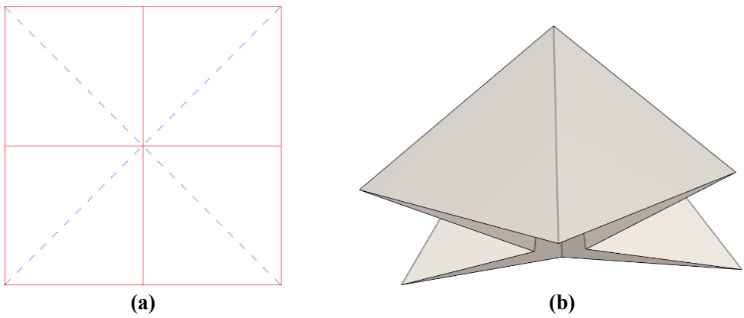

**Figure 1:** *(a) Crease pattern of a waterbomb base. (b) Folded waterbomb base*

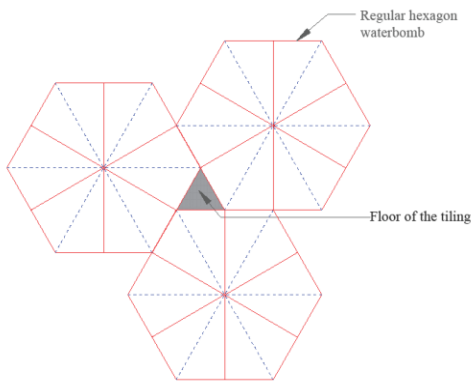

**Figure 2:** *Crease pattern of a regular hexagonal waterbomb tiling depicting the waterbomb molecule and the floor*

Origami tessellations such as the Miura-ori have been observed in historical fabric pleating dating back hundreds of years, see Figure 3 [Giegher 39], and it is possible that waterbomb tessellations were also considered in fabric, although we are unaware of any documented examples. The earliest example we can find dates back to the works of artist, computer scientist, and applied geometrist, Ron Resch [Resch 73]. Resch had asked himself, "What happens to a sheet as you wad it up?" and his explorations since 1961 appears to have led to his discovering the equilateral triangle waterbomb tessellation in this way. He then extrapolated those principles to a find square and a regular hexagon waterbomb tessellations.

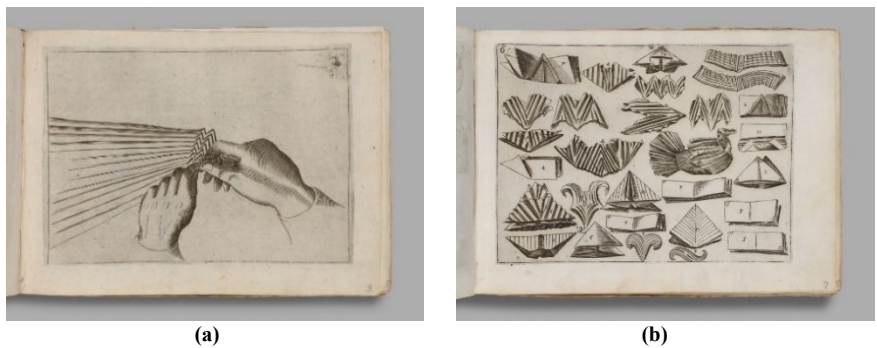

(a)  (b)

**Figure 3:** *Napkin folding from Li Tre Trattai by Mattia Giegher, 1639, showing how to fold the Miura-ori tessellation pattern in (a), and showing a waterbomb base in (b).*

Origami artists Katrin and Yuri Shumakov in 1996 tiled the square waterbomb but utilising only 2 hinge creases to form the model popularly known as the 'Magic Ball Pattern' [Shumakov and Shumakov 15].

More recently, waterbomb tessellations have been explored by artists such as Eric Gjerde, Benjamin Parker and Mitya Miller [Gjerde 09**,** Parker 11a, Miller 24], and Gjerde may be the first person to name the Waterbomb Tessellation in his book Origami Tessellations: Awe-Inspiring Geometric Designs [Gjerde 09].

We note that the Origamizer algorithm [Tachi 10, Demaine and Tachi 17] is able to produce CPs for a given pattern of "floors", which are not constrained to be co-planar as in our case. Some floor patterns can produce waterbomb molecules in Origamizer, but in general the algorithm produces creases that fall outside of waterbomb tessellations. And further, we aim to build up crease patterns starting from the molecules, while Origamizer works in the opposite direction, building the crease pattern from the floor pattern.

This paper aims to understand the properties of tilings with triangular waterbombs and extrapolate those learning to generate new patterns. We start by defining the restricted problem we wish to study and find criteria that the waterbomb tessellations must have. We then show many examples satisfying the criteria, including general infinite families for many of the patterns. We do not aim for a complete classification of such patterns in this paper. We then show examples of CPs having the same waterbomb molecule but different types of hinges active, and we end with a discussion about generalisations.

## 2   Waterbomb tessellation crease pattern criteria

To begin, we first define the most generalized waterbomb molecule to be any polygon with axial, ridge, and hinge creases given by Robert Lang's universal molecule algorithm [Lang 12]. In Figure 4 we define axial, ridge and hinge creases for a sample regular hexagon. The polygon's border creases are called axial creases. In the universal molecule algorithm, the polygon is in-set and shrunken with the vertices following along the angle bisectors of the original polygon, until at least two vertices coincide. At this point, ridge creases are drawn connecting the in-set vertices to the original vertices; that is, the ridge creases will be angle-bisectors. If there are more than one in-set vertices remaining, they are connected to each other, and the process is repeated until the new interior polygon has reduced to either a point or a line.

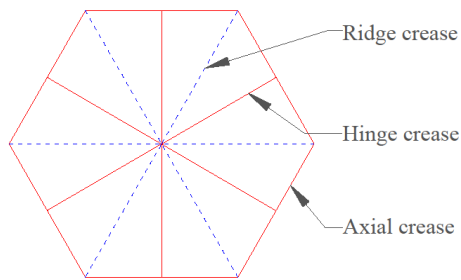

**Figure 4:** *Crease pattern of a regular hexagon waterbomb depicting ridge, hinge, and axial creases.*

From the ridge intersections, lines are drawn outward such that they intersect the polygon at right angles – these are the hinge creases. The hinge creases need not all be activated, and of those that are activated, they need not fold fully so that the axial crease on either side touch each other after folding. Instead, the hinge crease dihedral angles can be adjusted to suit the final folding of the waterbomb tessellation.

We are restricting our waterbomb molecules to be convex. As a further restriction of this paper, we are only interested in CPs tiled using only one waterbomb molecule, up to rotations and reflections and different hinge crease activations, although we briefly discuss the case of multiple molecules at the end. In addition, we require that after folding, the dihedral angle of the axial creases is , that is, the paper rises up at the floor edges perpendicularly.

If a waterbomb molecule tiles the entire CP without gaps, then the floors will be degenerate and consist only of the axial creases themselves. Otherwise, gaps in the CP between the waterbomb molecules will form polygonal floors with axial creases as their border.

The reason that we make use of the universal molecule algorithm is to ensure that the axial creases will all be co-planar after folding [Lang 12]. As a consequence, since the axial creases also form the edges of the floor polygons, all of the floors will also become co-planar after folding. In addition, since hinge creases fold axial creases back onto each other and the axial creases are co-planar after folding, the hinge creases must always be perpendicular to the axial creases, as is the case in the universal molecule. Likewise, the ridge creases must also fold axial creases onto each other, and so they must bisect the angles of the molecule, as in the universal molecule algorithm. Thus, the universal molecule is suitable for our construction of generalised waterbombs.

As a consequence of the universal molecule-defined waterbomb molecule creases, together with having floors in the CP, we have a few criteria that waterbomb tessellation CPs must satisfy:

1. In a waterbomb molecule adjacent to a floor, hinge creases must meet the axial creases at the vertices of the floor polygons.
2. In a waterbomb molecule adjacent to a floor, ridge creases must meet the axial creases at the vertices of the floor polygons.
3. Floors cannot share a border.
4. Floors can only meet at a point if a hinge crease also meets at that point and at least one ridge crease.
5. Floors must be convex.
6. Waterbomb molecules with $n$ sides must have $n$, $n$-1, or $n$-2 activated hinges.

We do not suppose that these are the only criteria that must be satisfied by waterbomb tessellations, but they can be useful in searching for new patterns.

To show the first result, if the hinge creases met the axial creases in the interior of an axial crease, it would cause the axial crease to be folded back onto itself, which in turn would force the adjacent floor polygon to also fold, which contradicts the definition of a floor as an unfolded polygon. The second result is true for the same reason.

Regarding result 3, if floors shared a border, then one of the floors would be folded with respect to the other one, and therefore the two floors would no longer be co-planar as required.

Regarding result 4, consider a circle immediately surrounding a vertex as in Figure 5 (a), where we suppose that there are two floors (shaded areas) touching at a point, with intervening axial and ridge creases. After folding the ridge creases, the two sets of axial creases will have joined their corresponding partners on the edge of one of the floors, with the result that if the floors remain in place, the paper will have ripped, as shown in Figure 5 (b). Therefore, two floors cannot meet at a vertex if there are only axial and ridge creases. In the case where two floors meet at a vertex where there is a hinge crease, then there can only be one hinge crease, since axial creases are at 90º with respect to each other when there is a hinge crease between them, and if there are two hinge creases, there would be no room for floors. Therefore, the remaining case is when there is one hinge crease, with necessarily at least one set of axial and ridge creases separating the two floors. The axial creases with a ridge bisector will bring the two floors to each other, and the rest of the paper will rise up at a right dihedral angle with the hinge connecting the two sides.

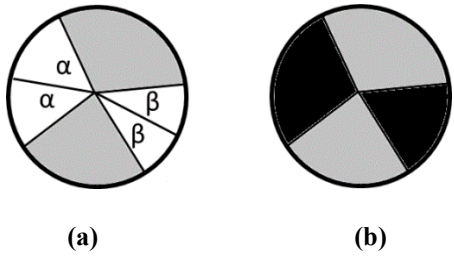

(a)  (b)

**Figure 5**: *Region of paper surrounding a vertex with two proposed floors (shaded) between axial and ridge creases (a) and the result of folding along the ridge creases (b) with ripped gaps shown in black.*

To see result 5, consider Figure 6, where there is a proposed concave floor, with a ridge crease bisecting the two axial creases meeting at the vertex. After folding, the paper must rise up perpendicularly to the floor, requiring that the angle θ be 90º. Since there are two angles θ, and the floor has an obtuse angle, the total angle around the vertex adds up to > 360º, which is not allowed in flat paper. There could be further axial and ridge creases meeting at the vertex beyond what is shown in Figure 6, however there must always be at least two faces rising up perpendicular to the floor, requiring two angles of 90º. Therefore, the floor angle cannot be obtuse, and so the floor polygons cannot be concave.

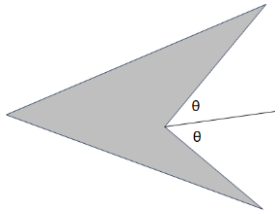

**Figure 6**: *A sample concave floor, shaded grey, with a ridge crease bisecting the angle between the two axial creases.*

To see result 6, consider Figure 7 showing the circular region around a folded waterbomb molecule (it could be of any number of sides, but a triangular one is shown). If the molecule is not equilateral, then consider a circle that contains all of the sides and assume that the faces of each flap are folded flush to each other. Seen from above, we see that the axial creases cross the circle so as to divide it into sectors. If one of the hinge creases is not activated, then an axial crease must remain unfolded and so it will form a diameter, creating a sector of 180º. At most, we can have two such sectors with their axial creases meeting along the diameter, and with the intervening flaps sandwiched between the diameter. Therefore, in a waterbomb molecule with *n* sides, either all hinges are activated, or else one or two may be left inactivated.

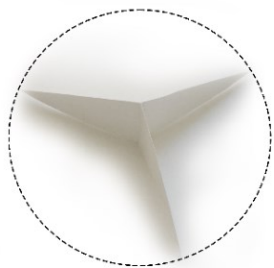

**Figure 7:** *View of a folded triangular waterbomb molecule from above.*

## 3 Patterns old and new

We now consider a variety of waterbomb tessellation CPs, all of which satisfy the criteria above.

### 3.1 Urchin Patterns

While most of the patterns below fold nicely and have only a small number of layers in each area of the folded model, there are a few very degenerate patterns, that we call urchin patterns. The name comes from the square urchin pattern whose flaps can be thinned to produce the various sea urchin origami models that are known, see e.g. [Montroll and Lang 91]. These urchin patterns have no floors except the axial creases themselves, but they are further degenerate in that all of the molecules overlap after folding, so that the number of overlapping folded layers grows as the pattern becomes larger. In addition, folding them requires some of the waterbomb molecules to have creases which are reversed compared to the rest. We show below in Figures 8-10 the three cases for triangle, square, and hexagon waterbomb molecules.

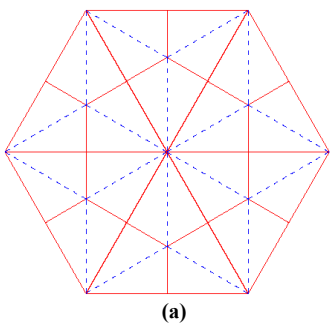
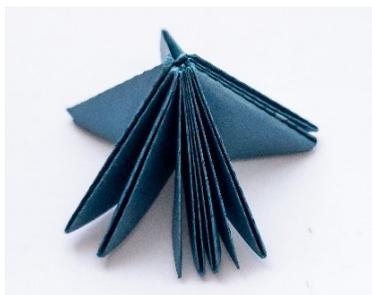

(a)          (b)

**Figure 8:** *(a) Crease pattern of a triangular urchin (b) Photo of the folded crease pattern*

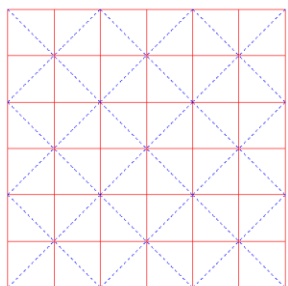
(a)

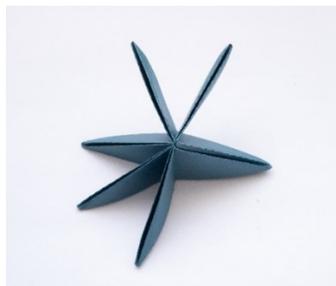
(b)

**Figure 9:** *(a) Crease pattern of a square urchin (b) Photo of the folded crease pattern*

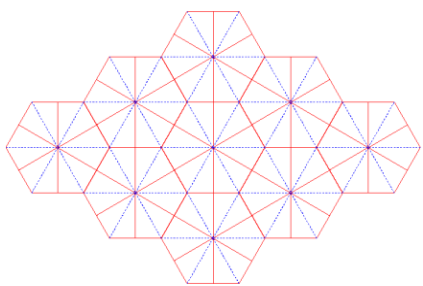
(a)

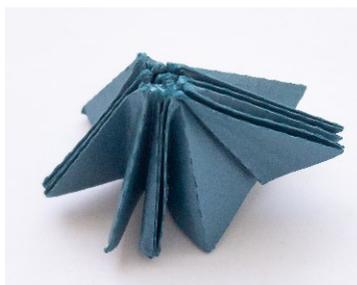
(b)

**Figure 10:** *(a) Crease pattern of a hexagonal urchin (b) Photo of the folded crease pattern*

We now move onto the main patterns of interest. The naming convention used is given by a one or two letter abbreviation of the molecule, followed by the number of hinge crease(s) activated in molecule, and then followed by the variation number. For example, T.2.1 denotes a Triangular waterbomb tessellation with two hinge creases activated and the first CP variation of this category.

### 3.2   Triangle:

Although Resch specifically constrained himself to equilateral triangles, see [Resch and Armstrong 70], we consider general scalene triangles in order to find novel tessellations.

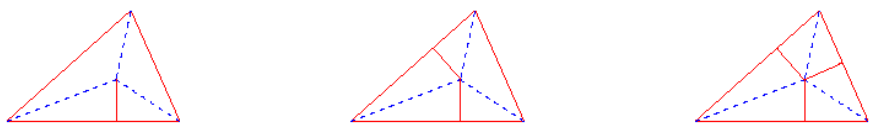

**Figure 11:** *Valid crease patterns of general triangle waterbomb showing ridge creases and 1, 2, and 3 activated hinge creases.*

### 3.2.1 With 1 hinge crease:

Pattern T.1: This pattern works for any triangle, see Figure 12. The tessellation does not have any floors and it has an overall angular tilt in general owing to the arrangement of ridge and hinge creases. However, in the case of isosceles and equilateral triangles (Figure 12 (c)-(d)), there is no tilt as the triangles are mirror symmetric. The equilateral triangle version of this pattern was previously discovered by Parker, see [Parker 11a].

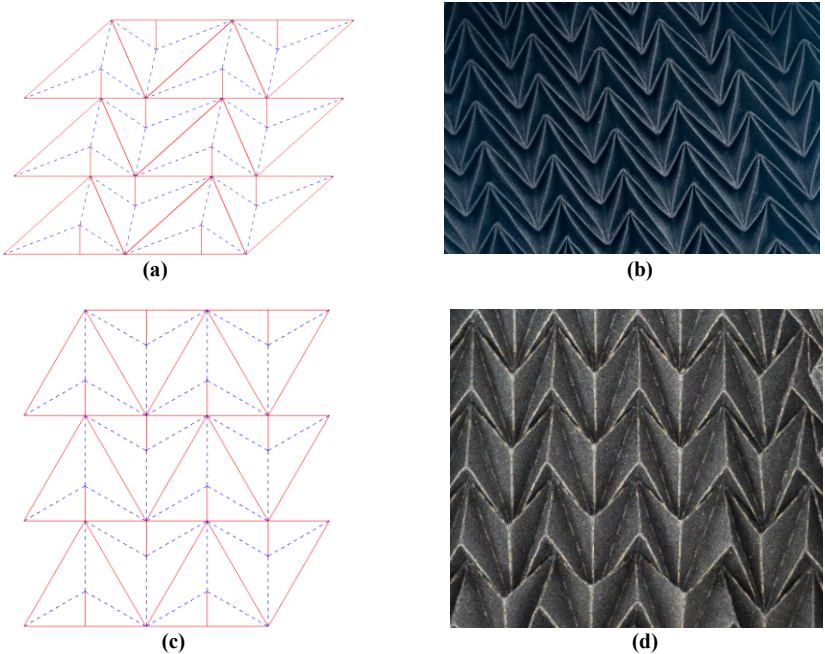

**Figure 12:** *Crease pattern of Pattern T.1 with scalene triangle waterbombs with one hinge crease (a) and photo of the folded crease pattern (b). Crease pattern of Pattern T.1 with equilateral triangle waterbombs with one hinge crease (c) and photo of the folded crease pattern (d).*

### 3.2.2 With 2 hinge creases:

Pattern T.2.1: This tessellation has parallelogram floors, see Figure 13. In the case of equilateral triangles, the floors become a rhombus, see Figure 13 (c) and (d). The equilateral triangle version of this pattern was previously discovered by Parker, see [Parker 11a].

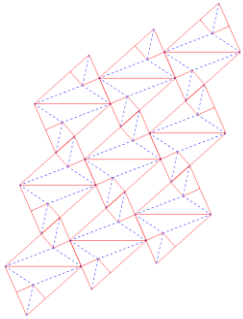
(a)

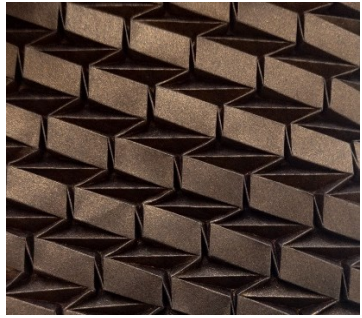
(b)

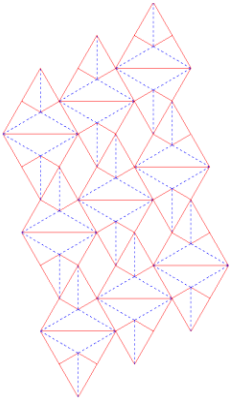
(c)

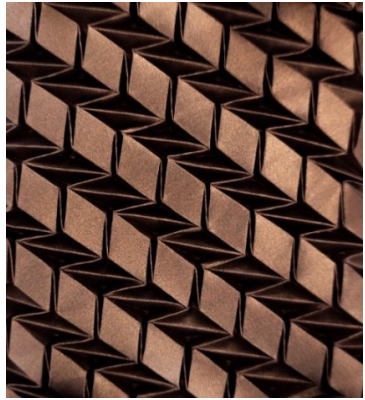
(d)

**Figure 13:** *Crease pattern of Pattern T.2.1 with scalene triangle waterbombs with two hinge creases (a) and photo of the folded crease pattern (b). Crease pattern of Pattern T.2.1 with equilateral triangle waterbombs with two hinge creases (c) and photo of the folded crease pattern (d).*

Pattern T.2.2: This tessellation has parallelogram floors, see Figure 14. It appears that this pattern has not been previously discovered.

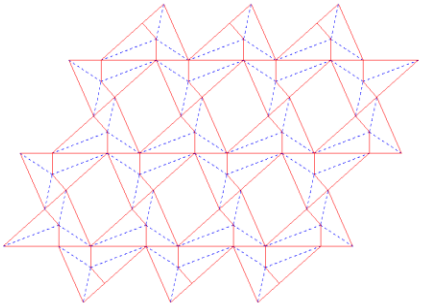
(a)

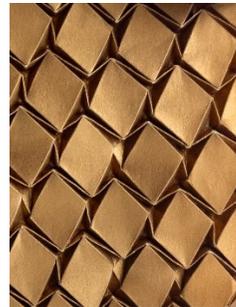
(b)

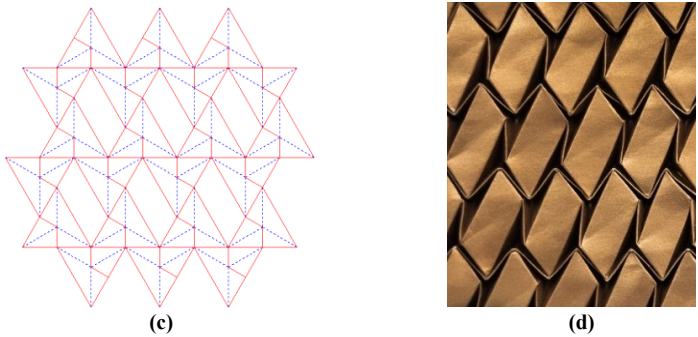

(c)           (d)

**Figure 14:** *Crease pattern of Pattern T.2.2 with scalene triangle waterbombs with two hinge creases (a) and photo of the folded crease pattern (b). Crease pattern of Pattern T.2.2 with equilateral triangle waterbombs with two hinge creases (c) and photo of the folded crease pattern (d).*

### 3.2.3    With 3 hinge creases:

Pattern T.3.1: This pattern generates hexagonal floors, the length of the sides of which is determined by the length between the vertices and the hinge creases, see Figure 15. Utilizing equilateral triangle waterbombs results in regular hexagonal floors, see Figure 15 (c)-(d). The equilateral triangle version of this pattern was previously discovered by Resch, see [Resch and Armstrong 70].

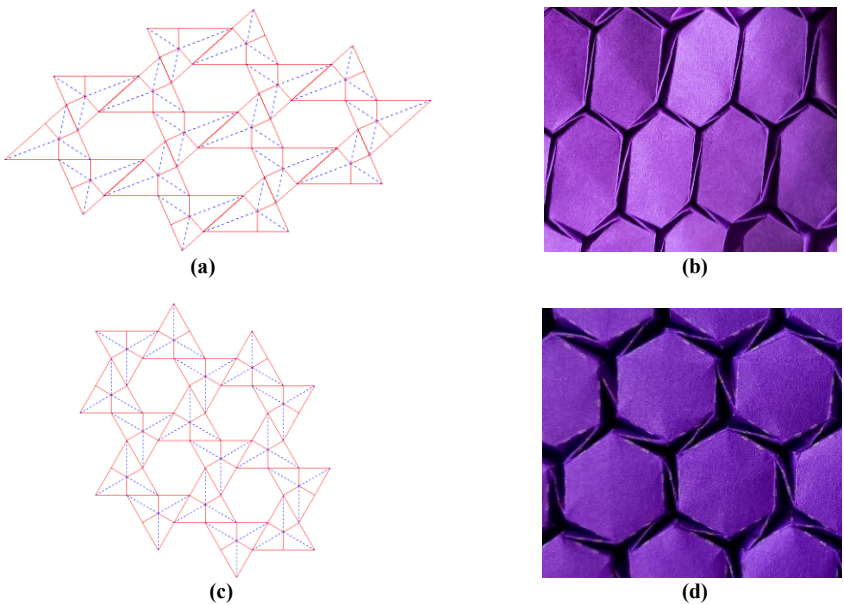

**Figure 15:** *Crease pattern of Pattern T.3.1 with scalene triangle waterbombs with three hinge creases (a) and photo of the folded crease pattern (b). Crease pattern of Pattern T.3.1*

with equilateral triangle waterbombs with three hinge creases (c) and photo of the folded crease pattern (d).

Pattern T.3.2: This tessellation has two different parallelogram floors for general triangles, see Figure 16 (a)-(b). However, in case of isosceles triangles the parallelogram floors are mirror symmetric, while for equilateral triangles, the floors are rhombuses, see Figure 16 (c)-(d). It appears that this pattern has not been previously discovered.

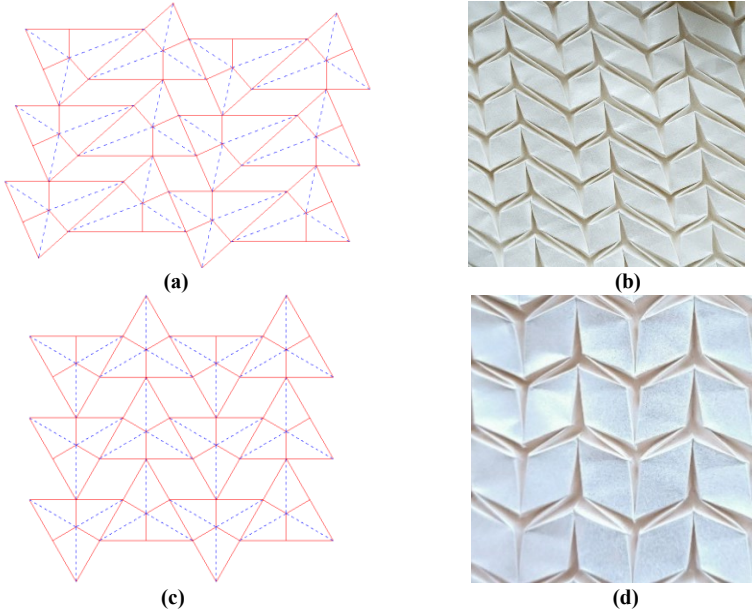

**Figure 16:** *Crease pattern of Pattern T.3.2 with scalene triangle waterbombs with three hinge creases (a) and photo of the folded crease pattern (b). Crease pattern of Pattern T.3.2 with equilateral triangle waterbombs with three hinge creases (c) and photo of the folded crease pattern (d).*

### 3.2.4 Equilateral triangle special cases:

We now consider waterbomb tessellations which require equilateral triangles and which cannot be generalised to other kinds of triangles.

Pattern E.T.3.1: This tessellation has floors which are equilateral triangles of half the size of the original unit, see Figure 17. Further, these floors are attached to their neighbouring floors only at each vertex. Resch discusses the creation of this tiling as the starting point for other patterns developed by him, utilising the same process. [Resch and Armstrong 70, Resch 73].

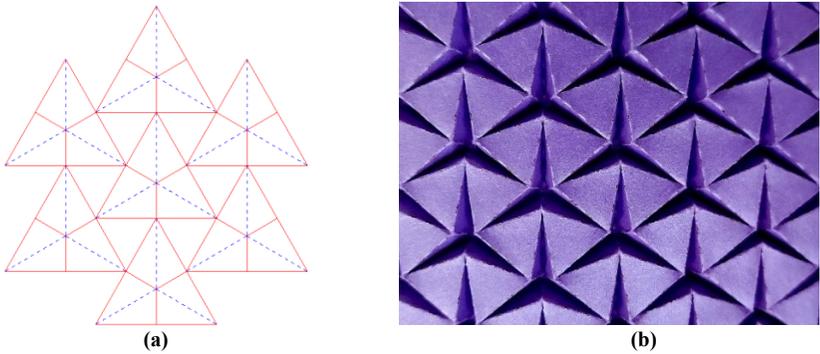

**Figure 17:** *Crease pattern of Pattern E.T.3.1 with equilateral triangle waterbombs with three hinge creases (a) and photo of the folded crease pattern (b)*

Pattern E.T.3.2: This tessellation has hexagonal floors with alternating equilateral triangles and rhombuses arranged around the periphery of hexagons, see Figure 18. This pattern was previously discovered by Parker, see [Parker 11a].

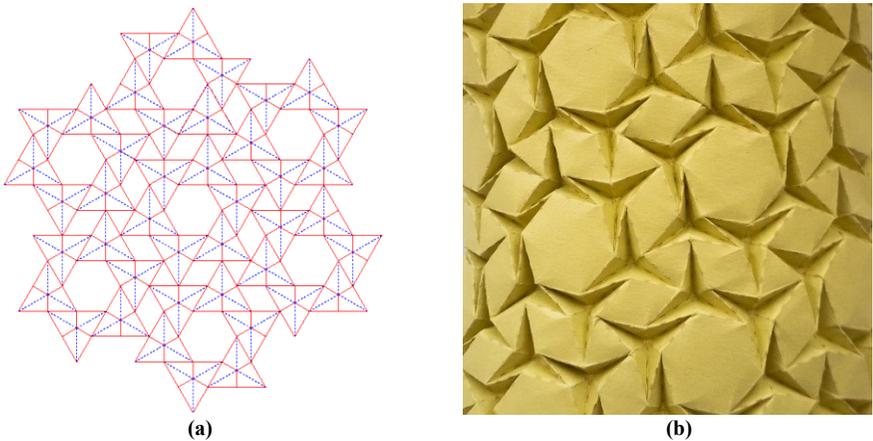

**Figure 18:** *Crease pattern of Pattern E.T.3.2 with equilateral triangle waterbombs with three hinge creases (a) and photo of the folded crease pattern (b).*

## 4 Patterns with multiple hinge activations

It is possible to combine in one pattern the same waterbomb molecule but with different hinge activations. We explore some such patterns in this section.

Pattern T.3.T.2: The molecule of Pattern T.3.1 can be shifted and tiled. This generates parallelogram floors along with the original hexagonal floors, see Figure 19. However, on doing so, there are positions where the hinge creases do not line up, as depicted in the incorrect CP of Figure 19. So, in positions where the hinge creases do not join to make a valid tiling, we can remove them to create a valid tiling which has parallelogram floors. Depending on the way the molecules are tiled, different symmetries with varied orientation of the parallelogram floors can be obtained, see Figure 20.

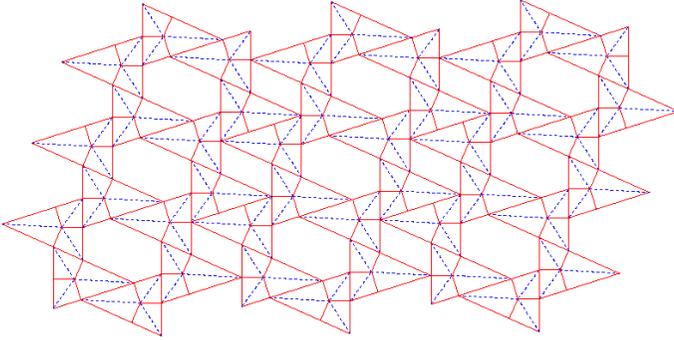

**Figure 19:** *Crease pattern of a possible tiling with invalid hinge creases at certain positions*

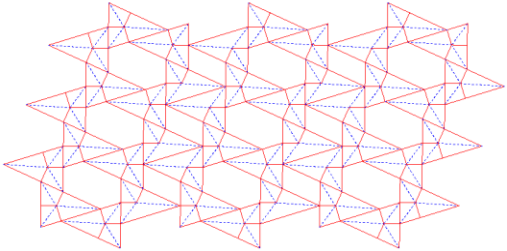

(a)

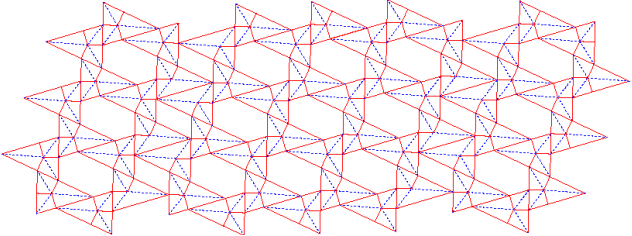

(b)

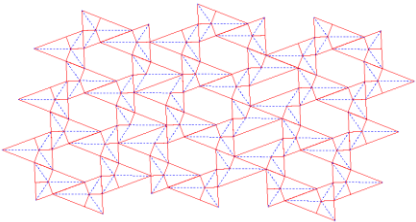

(c)

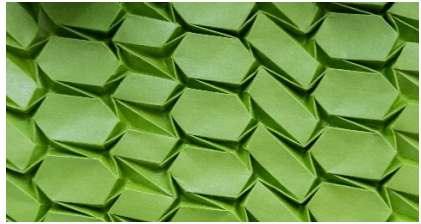

(d)

**Figure 20:** *Crease patterns of possible tiling with corrected hinge creases at certain positions (a), (b) and (c) and photo of folded crease pattern of (c) with gaps separating the floors (d).*

These examples bring up another issue, in that after folding, gaps may appear that separate the floors. A full exploration of these patterns requires further work.

## 6 Discussion

We have compiled above many triangular waterbomb tessellation CPs with our restriction of using only one waterbomb molecule. Most of the patterns we found involve entire families of CPs. We believe that some of these patterns are new, and even in the cases where they were previously known, our scalene versions may be new. In addition, we have only briefly considered patterns having one waterbomb molecule but with multiple hinge activations. Future work is needed to explore these further. The present compilation of examples is certainly missing patterns.

We note when looking at the patterns above that for the CPs featuring waterbomb molecules with 2 hinges, the folded patterns can feature extra creases between the floors, connecting them, which adds bulk to the collapse between the floors. We also note that pattern with 1 hinge was degenerate, without floors. We are not sure if these are features of waterbomb tessellations in general.

We saw above that the only triangle pattern which is restricted to work only for equilateral triangles coincides with the only pattern where none of the waterbomb molecule edges lie adjacent to another waterbomb. This does not appear to generalise, since for the regular hexagon molecule, the CP requires molecules to share borders.

We believe that the current work is just the start of exploring waterbomb tessellations more systematically, and we hope that the criteria found above can help to explore this space further.

## References

[Demaine and Tachi 17] Erik Demaine and Tomohiro Tachi. 2017. "Origamizer: A practical algorithm for folding any polyhedron." 10.4230/LIPIcs.SoCG.2017.34. https://hdl.handle.net/1721.1/137685

[Giegher 39] Mattia Giegher. "Li Tre Trattati Di Messer". P. Frambotto, Padova, Italy. 1639.

[Gjerde 09] Eric Gjerde. "Origami Tessellations: Awe-Inspiring Geometric Designs". A K Peters, Wellesley, MA. 2009.

[Lang 12] Robert J. Lang. "Origami Design Secrets: mathematical methods for an ancient art - 2nd ed.". CRC Press, Boca Raton, FL. 2012.

[Miller 24] Mitya Miller. "mityamiller". 2024. https://www.instagram.com/mityamiller/

[Montroll and Lang 91] John Montroll and Robert J. Lang. "Origami Sea Life". Dover Publications Inc., NY.1991.

[Parker 11a] Benjamin Parker. "Triangular WBT's". 03 November, 2011. https://www.flickr.com/photos/brdparker/6307774095/in/album-72157603945955019/

[Parker 11b] Benjamin Parker. "Hexagonal WBT's". 03 November, 2011. https://www.flickr.com/photos/brdparker/6308295602/in/album-72157603945955019/

[Parker 11c] Benjamin Parker. "Square WBT's". 03 November, 2011. https://www.flickr.com/photos/brdparker/6308295590/in/album-72157603945955019/

[Parker 11d] Benjamin Parker. "Rectangular WBT's". 03 November, 2011. https://www.flickr.com/photos/brdparker/6307774025/in/album-72157603945955019/

[Resch and Armstrong 70] Ronald D. Resch and Elmer Armstrong. "The Ron Resch Paper and Stick Film. Resch Films", 1970

[Resch 73] , Ronald D. Resch. "The topological design of sculptural and architectural systems." Proceedings of the June 4-8, 1973, national computer conference and exposition. 1973

[Shumakov and Shumakov 15] Yuri Shumakov and Katrin Shumakov. "Origami Magic Ball Wonders: From Dragon's Egg to Hot Air Balloon". CreateSpace Independent Publishing Platform, 2015

[Tachi 10] Tomohiro Tachi. "Origamizing polyhedral surfaces." IEEE Transactions on Visualization and Computer Graphics, 16(2):298–311, 2010. doi:10.1109/TVCG.2009.67.